# Sectoral and spatial decomposition methods for multi-sector capacity expansion models

*Federico Parolin [a,b], Yu Weng [b], Paolo Colbertaldo [a], Ruaridh Macdonald [b]*

a. Department of Energy, Politecnico di Milano, Milan, Italy
b. MIT Energy Initiative, Massachusetts Institute of Technology, Cambridge, MA, USA

## Abstract

Multi-sector capacity expansion models play a crucial role in energy planning by providing decision support for policymaking in technology development. To ensure reliable support, these models require high technological, spatial, and temporal resolution, leading to large-scale linear programming problems that are often computationally intractable. To address this challenge, conventional approaches rely on simplifying abstractions that trade accuracy for computational efficiency. Benders decomposition has been widely explored to improve computational efficiency in electricity capacity expansion models. Specifically, state-of-the-art methods have primarily focused on improving performance through temporal decomposition. However, multi-sector models introduce additional complexity, requiring new decomposition strategies. In this work, we propose a budget-based formulation to extend decomposition to the sectoral and spatial domains. We test the developed sectoral and spatial Benders decomposition algorithms on case studies of the continental United States, considering different configurations in terms of spatial and temporal resolution. Results show that our algorithms achieve substantial performance improvement compared to existing decomposition algorithms, with runtime reductions within 15%-70%. The proposed methods leverage the generic structure of multi-sector capacity expansion models, and can thus be applied to most existing energy planning models, ensuring computational tractability without sacrificing resolution.

## Highlights

- We apply Benders decomposition to multi-sector capacity expansion models
- We develop sectoral and spatial Benders decomposition algorithms
- The algorithms achieve 15%-70% faster runtimes that existing decomposition methods

## 1. Introduction

Capacity expansion models (CEMs) are essential tools to support policymakers in the design of technically feasible and cost-effective decarbonisation strategies (Savvidis et al., 2019). Most CEMs are formulated as deterministic linear programming (LP) problems, with runtime and memory usage scaling quadratically with model size (Ringkjøb et al., 2018). The urgency of climate action and the need for economy-wide decarbonisation have driven researchers to develop increasingly complex capacity expansion models, incorporating multiple sectors (Chang et al., 2021), increasing the spatial resolution (Xiong et al., 2024), and extending the time horizon to several years (Bhatt, 2025; Ruggles et al., 2024). As a result, modern capacity expansion models risk computational intractability, requiring the introduction of substantial abstractions (A. F. Jacobson et al., 2024).

One of the most common approaches to limit computational complexity is temporal aggregation (Hoffmann et al., 2020), which consists of reducing the number of modelled time steps by aggregating input time series into representative periods (Prina et al., 2020). While temporal aggregation enables significant reduction of the computational load (Kotzur et al., 2018), it has been shown to introduce non-negligible errors (Fleschutz et al., 2025). Nonetheless, many studies heavily rely on temporal aggregation (Ringkjøb et al., 2018). For example, Larson et al. (2021) model a 35-year transition of the United States using 8 snapshots, each modelled with 28 representative days, covering 2% of the total time horizon. Similarly, the clustering of geographic areas into large zones is frequently used to limit the problem size and complexity (Hoffmann et al., 2024; Javanmardi et al., 2025). However, high spatial resolution has been shown to be critical for identifying transmission bottlenecks and accounting for variable renewable availability (Sahin et al., 2024). While these simplifications improve computational tractability, they can significantly impact model accuracy, limiting the validity and generalisability of results (A. F. Jacobson et al., 2024). In addition, CEMs may remain intractable even with the implementation of such simplifications, due to the large scale of the problems they address.

The use of decomposition methods is a common alternative approach to enable solving of complex models, reducing the need to introduce abstractions. Among these, Benders Decomposition (BD) is a widely used approach that separates investment decisions from operational dispatch decisions (Benders, 1962). As schematised in Figure 1, BD iteratively solves two optimisation problems: an upper problem which determines investment decisions, and a lower problem which optimises operational decisions over a planning horizon, fixing investment decisions as resulting from the upper problem solution. The lower problem is solved to generate linear inequalities (i.e., constraints) typically referred to as "cuts", which are added to the upper problem. Typically, the lower problem represents a year of operation with hourly resolution.

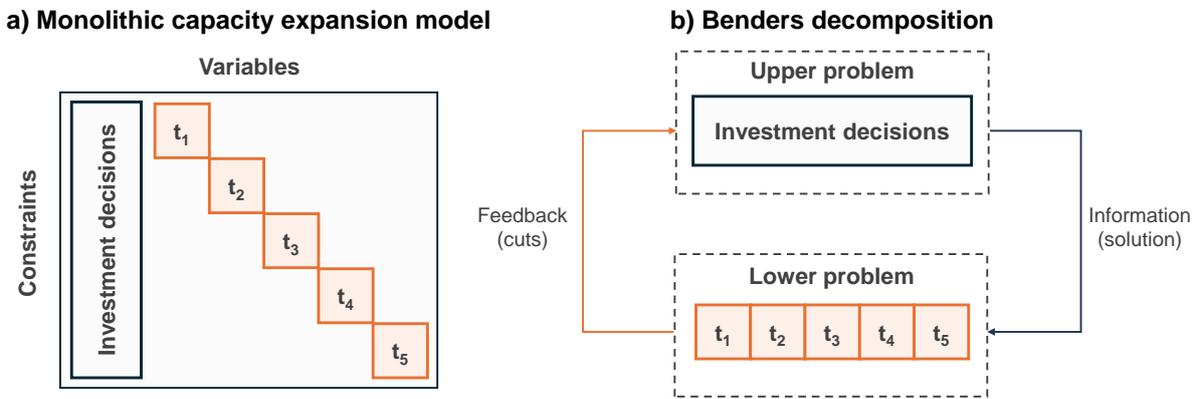

Figure 1. Generic block structure of a monolithic capacity expansion model (a) and schematised structure of Benders decomposition (b). Orange blocks corresponds to operational decisions for a zone/generator in the relative time step.

Recent studies have explored various refinements and applications of BD to improve the computational efficiency of energy system models. Lara et al. (2018) applied BD to multi-period electric infrastructure planning, decomposing the planning horizon into time periods, so that the lower problem is further divided into a set of subproblems, each corresponding to a planning period. The presence of multiple subproblems instead of a single lower problem is generally beneficial for convergence, as their solution enables the generation of multiple cuts at each iteration. In the work by Lara et al., however, subproblems represent an entire year of operation, still requiring temporal aggregation to manage the computational complexity. In their application, temporal coverage is limited to 12 typical days per year. Li et al. (2022) further refined this approach by aggregating complicating variables into a single upper problem, enabling independent solution of subproblems. Nevertheless, tractability issues persist for subproblems, requiring temporal aggregation with a maximum resolution of 15 typical days per year. Munoz et al. (2016) developed a novel approach that enables the parallel solution of subproblems in BD. However, their approach requires long-duration storage and policy constraints to be neglected, as subproblems are independent in time. Mazzi et al. (2021) also proposed an approach that solves subproblems in parallel. Within the context of a stochastic investment planning problem for power systems, they select a subset of subproblems to be solved exactly, while the remainder are approximated using oracles. Göke et al. (2024) introduced a regularisation (or stabilisation) step to improve the computational efficiency of BD, considering applications to two-stage stochastic problem for power grid planning. Jacobson et al. (2024) proposed a novel BD-based method for deterministic CEMSs that further decomposes the problem temporally, dividing a single-year time horizon into shorter subperiods, resulting in a set of subproblems that can be solved in parallel. The presence of multiple subproblems resulting from this "temporal" decomposition enables the generation of multiple cuts per iteration, improving convergence. In addition, the relatively small size of subproblems allows for their parallel solving, reducing the runtime per iteration. Compared to other similar approaches (Munoz et al., 2016), specific constraints are introduced to model time-coupled features such as long-duration storage and emission constraints. Pecci and Jenkins (2024) extended the approach to the case of multi-period models and introduced a regularisation step, further improving convergence.

Despite these advancements, most applications of BD remain focused exclusively on the power sector, overlooking the complexities of integrated energy systems. Even within single-sector applications, Jacobson et al. (2024) demonstrated that computational time scales quadratically with the number of spatial zones considered. Hence, multi-sector CEMs with high spatial resolution may still face significant computational constraints. To address this challenge, we extend the temporal decomposition approach proposed by Jacobson et al. (2024) and Pecci and Jenkins (2024) to multi-sector CEMs. As temporal decomposition alone might not adequately address the increased computational complexity of multi-sector models, we introduce a novel budget-based formulation to efficiently decompose the problem, first sectorally, and then spatially. The developed method leverages the generic structure of capacity expansion models, and can thus be implemented to most existing energy system models. In this work, we implement the developed methods in the Julia-based Dolphyn model (He et al., n.d.).

The remainder of this work is structured as follows. Section 2 presents the generic mathematical formulation of integrated energy system models, introducing the nomenclature used throughout the chapter. In Section 3, the temporal BD is extended to multi-sector models, while Section 4 introduces the developed budget-based formulation for sectoral and spatial BD. The computational performance of the proposed decomposition algorithms are evaluated in Section 5 compared to monolithic formulations and the state-of-the-art temporal BD. Finally, Section 6 summarises the main findings and discusses potential directions for future developments enabled by the proposed methods.

## 2. Problem formulation

The generic formulation of multi-sector capacity expansion models is introduced here to highlight the underlying structure upon which the proposed decomposition methods are developed. Let the operational year be divided into $w \in W$ subperiods (e.g., days or weeks), the analysed geographic area into $z \in Z$ zones, and the integrated energy system into $s \in S$ sectors. By defining $y \in \mathbb{R}^m$ as the vector of all the investment decision variables, $x_{w,z,s} \in \mathbb{R}^n$ as the vector of all operational decision variables of sub-period $w$, zone $z$, and sector $s$, and $c_y^T$ and $c_{w,z,s}^T$ as the vectors of annualised investment and operational costs, respectively, the optimisation problem can be formulated as:

$$\min c_y^T y + \sum_{w \in W} \sum_{z \in Z} \sum_{s \in S} c_{w,z,s}^T x_{w,z,s} \tag{1a}$$

$$\text{s.t.} \quad A_{w,z,s} x_{w,z,s} + B_{w,z,s} y \leq b_{w,s,z} \quad \forall w \in W, z \in Z, s \in S \tag{1b}$$

$$\sum_{w \in W} \sum_{z \in Z} \sum_{s \in S} Q_{w,z,s} x_{w,z,s} \leq e \tag{1c}$$

$$Ry \leq r \tag{1d}$$

$$y \geq 0 \tag{1e}$$

$$x_{w,z,s} \geq 0 \quad \forall w \in W, z \in Z, s \in S \tag{1f}$$

where the matrices $A_{w,z,s}$ and $B_{w,z,s}$ and the vector $b_{w,s,z}$ are defined such that Eq. (1b) encompasses all balances and operational constraints, the matrix $Q_{w,z,s}$ and the vector $e$ are defined such that Eq. (1c) corresponds to the net emission constraint, and the matrix $R$ and the vector $r$ are defined such that Eq. (1d) represents all investment constraints.

## 3. Temporal Benders decomposition

We extend the approach proposed by Jacobson et al. (2024) and Pecci and Jenkins (2024), referred to as temporal decomposition hereafter, to the case of a multi-sector capacity expansion model. The temporal decomposition approach is graphically schematised in Figure 2. Building on conventional Benders decomposition, the operational problem is further decomposed in time, generating a subproblem for each subperiod $w \in W$. This structure enables the use of parallel computing, thus reducing the computational time per iteration. Temporal aggregation can be applied to limit the number of subperiods.

The budgeting variables $q_w$ are introduced as complicating variables[1] to model the net emission constraint, which would otherwise link all time steps and prevent separability. Similarly, the storage level and the change in storage level across each subperiod of long-duration storage technologies are treated as investment decisions and included in the complicating variables $y \in Y$. When temporal aggregation is applied, violations of state of charge limits in non-representative periods are prevented by implementing the formulation developed in Parolin et al. (2024).

---

[1] Complicating variables are decision variables that establish a dependence between the upper and lower problems, and are therefore fixed when solving subproblems. In the standard BD, they correspond to investment decision variables.

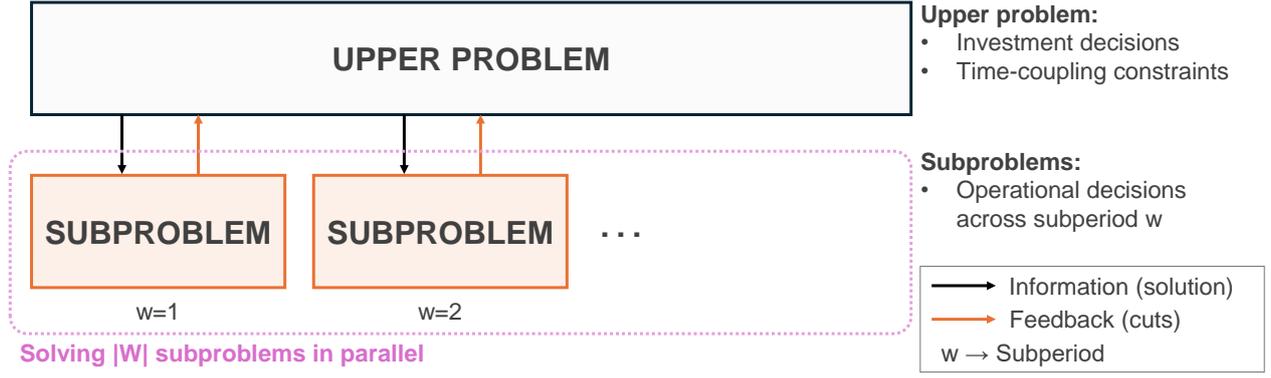

*Figure 2. Schematised structure of the temporal Benders decomposition* (A. Jacobson et al., 2024; Pecci & Jenkins, 2024).

At an iteration k, given a choice of investment decisions $y^k$ and budgeting variables $q_w^k$, the following operational subproblem is solved for each subperiod $w \in W$:

$$g_w^k = \min \sum_{z \in Z} \sum_{s \in S} c_{w,z,s}^T x_{w,z,s} \tag{2a}$$

$$\text{s.t.} \quad A_{w,z,s} x_{w,z,s} + B_{w,z,s} y \leq b_{w,s,z} \quad \forall z \in Z, s \in S \tag{2b}$$

$$\sum_{z \in Z} \sum_{s \in S} Q_{w,z,s} x_{w,z,s} \leq q_w \tag{2c}$$

$$y = y^k \quad : \pi^k \tag{2d}$$

$$q_w = q_w{}^k \quad : \lambda_w^k \tag{2e}$$

$$x_{w,z,s} \geq 0 \quad \forall z \in Z, s \in S \tag{2f}$$

where $\pi^k$ and $\lambda_{w,z,s}^k$ are the Lagrangian multipliers of the corresponding constraints. By solving subproblems (2), the best upper bound at iteration k can be computed:

$$UB^k = \min_{j=0,\ldots,k} c_y^T y^j + \sum_{w \in W} g_w^j \tag{3}$$

The lower bound is determined by solving the upper problem, which uses the Lagrangian multipliers obtained by solving subproblems to add a cut for each $w \in W$. The resulting problem is:

$$LB^k = \min c_y^T y^j + \sum_{w \in W} \theta_w \tag{4a}$$

$$\text{s.t.} \quad \theta_w \geq g_w^j + (y - y^j)^T \pi^j + (q_w - q_w^j)^T \lambda_w^j \quad \forall j = 0, \ldots, k, w \in W \tag{4b}$$

$$\sum_{w \in W} q_w \leq e \tag{4c}$$

$$Ry \leq r \tag{4d}$$

$$y \in Y \tag{4e}$$

The new estimates of investment decisions and budgeting variables are obtained by solving the following regularisation problem of Eq. (5). Among the different methods analysed in Pecci & Jenkins (2024), an interior point strategy is selected as it emerged as the best-performing option from preliminary simulations performed on multi-sector systems, considering a level-set parameter α of 0.5. The resulting algorithm is presented in Algorithm 1.

$$\min \Phi^{int} = 0 \tag{5a}$$

$$\text{s.t.} \quad \theta_w \geq g_w^j + (y - y^j)^T \pi^j + (q_w - q_w^j)^T \lambda_w^j \quad \forall j = 0, \ldots, k, w \in W \tag{5b}$$

$$\sum_{w \in W} q_w \leq e \tag{5c}$$

$$Ry \leq r \tag{5d}$$

$$y \in Y \tag{5e}$$

$$c_y^T y^j + \sum_{w \in W} \theta_w \leq LB^k + \alpha(UB^k - LB^k) \tag{5f}$$

---

*Algorithm 1.* Temporal Benders decomposition (A. Jacobson et al., 2024; Pecci & Jenkins, 2024).

**Input:** $y^0 = 0, q_w^0 = 0 \ \forall w \in W$. Set maximum number of iterations $K_{max}$ and convergence tolerance $\varepsilon_{tol}$.
**Output:** $y^{opt} = 0, q_w^{opt} = 0 \ \forall w \in W$
**for** $k = 0, \ldots, K_{max}$ **do**
    **for** $w \in W$ **do**
        Solve operational subproblem (2).
    **end for**
    Compute best upper bound $UB^k$ as in Eq. (3).
    Update cuts in upper problem as in Eq. (4b).
    Solve upper problem (4) to obtain lower bound $LB^k$.
    **if** $(UB^k - LB^k)/LB^k \leq \varepsilon_{tol}$ **then**
        Set $y^{opt} = y^{k+1}$ and $q_w^{opt} = q_w^{k+1} \ \forall w \in W$
        **stop**
    **else**
        Solve regularised upper problem (5) to obtain $y^{k+1}$ and $q_w^{k+1}$
    **end if**
**end for**

---

## 4. Budget-based formulation for sectoral and spatial Benders decomposition

The structure of multi-sector capacity expansion models can be exploited to further decompose the operational problem sectorally or spatially. Sectors and zones are generally loosely coupled, the only linking variable being the exchange of an energy vector, thus allowing for effective decomposition. As Figure 3 shows, this approach leads to the generation of one operational subproblems for each combination of subperiod and sector or subperiod and zone for sectoral and spatial decomposition, respectively. As a result, the computational time per iteration can be reduced by further exploiting distributed computing, solving smaller subproblems in parallel, and the total number of iterations can be decreased thanks to the generation of additional cuts (one per subproblem).

### a) Temporal + sectoral Benders decomposition

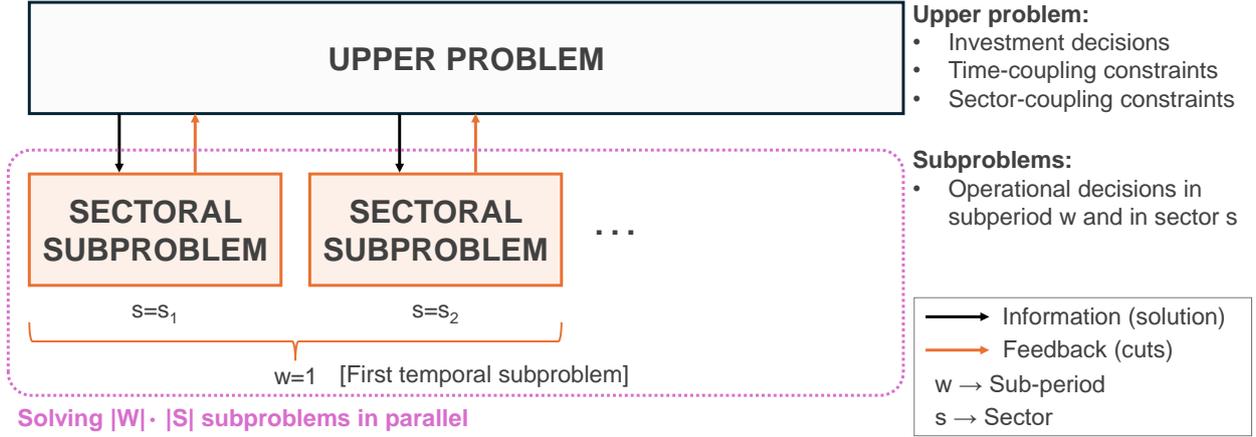

### b) Temporal + spatial Benders decomposition

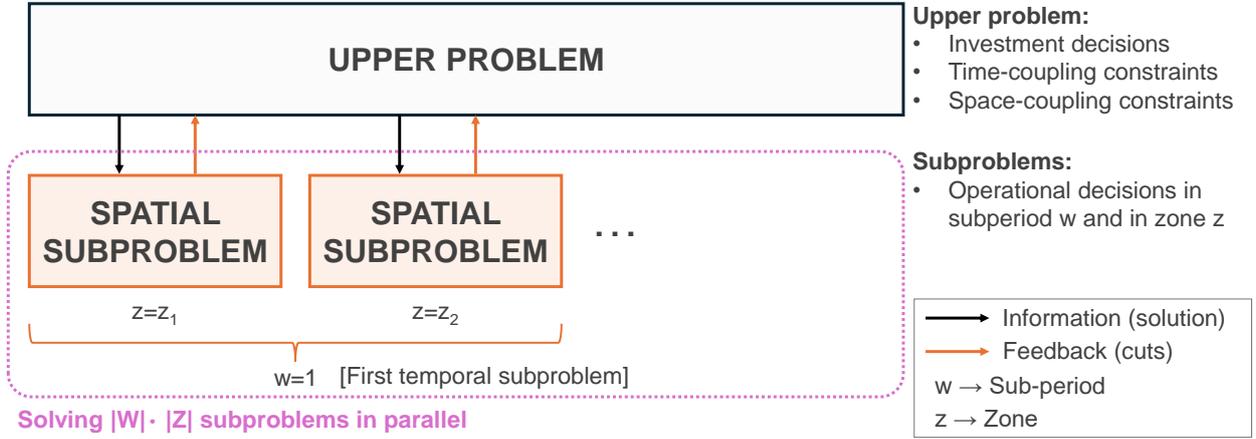

Figure 3. Schematised structure of the sectoral (a) and spatial (b) Benders decomposition.

The balance of an energy vector v in sector s at zone z and time step t can be expressed highlighting the role of energy vector exchanges with other sectors and zones, as:

$$\sum_g x_{gen}^{g,v,z,s,t} + \sum_\sigma \frac{x_{otp}^{\sigma,v,z,s,t}}{\eta_{\sigma,otp}} + x_{nse}^{v,z,s,t} - x_{dem}^{v,z,s,t} - \sum_\sigma x_{ipt}^{\sigma,v,z,s,t} \cdot \eta_{\sigma,ipt} - x_{crt}^{v,z,s,t} = \sum_{s' \neq s} x_{exp}^{v,z,s,s',t} + \sum_{z' \neq z} x_{trn}^{v,z,z',t} \qquad (6)$$

where $x_{gen}^{g,v,z,s,t}$ is the energy vector generation from technology g, $x_{otp}^{\sigma,v,z,s,t}$ and $x_{ipt}^{\sigma,v,z,s,t}$ are the output and input flow from storage technology $\sigma$, $x_{nse}^{v,z,s,t}$ is the non-served energy[2], $x_{dem}^{v,z,s,t}$ is the energy vector demand, $x_{crt}^{v,z,s,t}$ is the curtailment, $x_{exp}^{v,z,s,s',t}$ is the net energy vector export from sector s to sector s', and $x_{trn}^{v,z,z',t}$ is the transport flow from zone z to zone z'.

Sectoral export flows ($x_{exp}^{v,z,s,s',t}$) and transport flows ($x_{trn}^{v,z,z',t}$) are treated as complicating variables to decompose the problem sectorally and spatially, respectively. However, preliminary simulations showed that linking flows at hourly resolution significantly hinders convergence, failing to deliver the expected benefits of decomposition. We therefore introduce a budget-based linking, using budgets for every subperiod as complicating variables instead of hourly flows, as schematised in Figure 4. If $T_w$ is the index set of all time steps within subperiod w, export budgets for sectoral decomposition are defined as:

$$\forall z, w, s \neq s' \qquad y_{exp}^{v,z,s,s',w} = \sum_{t \in T_w} x_{exp}^{v,z,s,s',t} \qquad (7)$$

while transport budgets for spatial decomposition are defined as:

$$\forall z, v, w \qquad y_{trn}^{v,z,w} = \sum_{t \in T_w} \sum_{z' \neq z} x_{trn}^{v,z,z',t} \qquad (8)$$

---

[2] Non-served energy is defined as the amount of consumption that is not covered endogenously and is associated with a penalty cost. In Benders decomposition, it also acts as slack variable to make subproblems feasible in the absence of sufficient generation capacity.

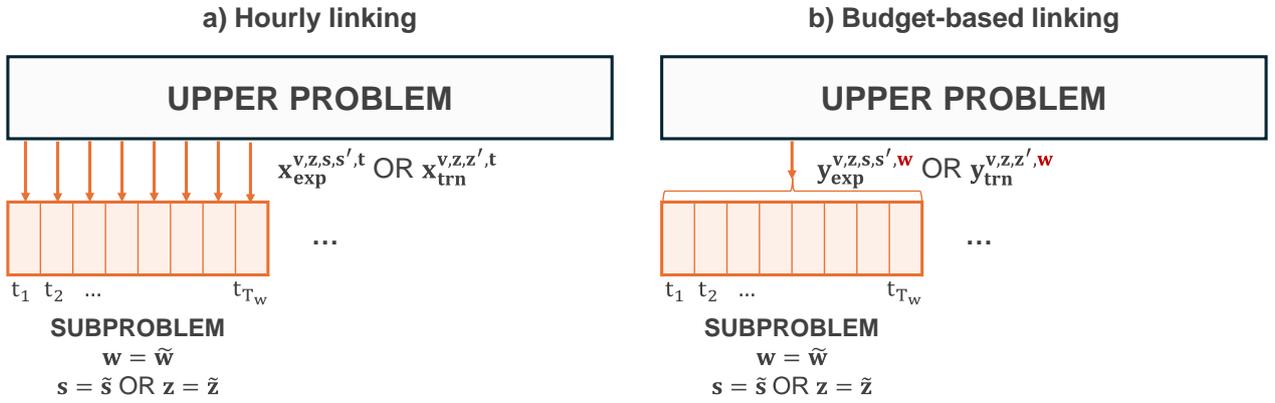

Figure 4. Hourly vs budget-based linking between upper problem and subproblems.

Budgets are treated as investment decision variables and included in the set $y \in Y$. Problem (2) is therefore modified as in problem (9), where $\delta = s$ and $\varphi = z$ in the sectoral decomposition and $\delta = z$ and $\varphi = s$ in the spatial decomposition.

$$g_{w,\delta}^k = \min \sum_{\varphi} c_{w,\varphi,\delta}^T x_{w,\varphi,\delta} \qquad (9a)$$

$$\text{s.t.} \quad A_{w,\varphi,\delta} x_{w,\varphi,\delta} + B_{w,\varphi,\delta} y \leq b_{w,\delta,\varphi} \quad \forall \varphi \qquad (9b)$$

$$\sum_{\varphi} Q_{w,\varphi,\delta} x_{w,\varphi,\delta} \leq q_{w,\delta} \qquad (9c)$$

$$y = y^k \quad : \pi^k \qquad (9d)$$

$$q_{w,s} = q_w^k \quad : \lambda_{w,\delta}^k \qquad (9e)$$

$$x_{w,\varphi,\delta} \geq 0 \quad \forall \varphi \qquad (9f)$$

In the sectoral decomposition, the operational subproblem (9) is solved for each subperiod and sector, adding Eq. (10) to ensure consistency between budgets and hourly flows. Similarly, in the spatial decomposition the operational subproblem (9) is solved for each subperiod and zone with the addition of Eq. (11).

$$\forall z, v, w \qquad \sum_{t \in T_w} x_{\exp}^{v,z,s',t} = y_{\exp}^{v,z,s',s,w} \quad \forall z \in Z \tag{10}$$

$$\forall z, v, w \qquad \sum_{t \in T_w} \sum_{z' \neq z} x_{trn}^{v,z,z',t} = y_{trn}^{v,z,w} \quad \forall v \in V \tag{11}$$

The new upper problem is defined as:

$$LB^k = \min c_y^T y^j + \sum_w \sum_\delta \theta_{w,\delta} \tag{12a}$$

$$\text{s.t.} \quad \theta_{w,\delta} \geq g_{w,\delta}^j + (y - y^j)^T \pi^j + \left(q_{w,\delta} - q_{w,\delta}^j\right)^T \lambda_{w,\delta}^j \quad \forall w, \delta, j = 0, \ldots, k \tag{12b}$$

$$\sum_w \sum_\delta q_{w,\delta} \leq e \tag{12c}$$

$$Ry \leq r \tag{12d}$$

$$y \in Y \tag{12e}$$

The regularisation problem is updated accordingly, as:

$$\min \Phi^{int} = 0 \tag{13a}$$

$$\text{s.t.} \quad \theta_{w,\delta} \geq g_{w,\delta}^j + (y - y^j)^T \pi^j + \left(q_{w,\delta} - q_{w,\delta}^j\right)^T \lambda_{w,\delta}^j \quad \forall w, \delta, j = 0, \ldots, k \tag{13b}$$

$$\sum_w \sum_\delta q_{w,\delta} \leq e \tag{13c}$$

$$Ry \leq r \tag{13f}$$

$$y \in Y \tag{13g}$$

$$c_y^T y^j + \sum_w \sum_\delta \theta_{w,\delta} \leq LB^k + \alpha(UB^k - LB^k) \tag{13h}$$

The solution algorithm is reported in Algorithm 2. In this section we introduced a generalised formulation for budget-based decomposition algorithms, which can be adapted for both sectoral and spatial BD. The detailed separate formulations of the two decomposition algorithms are reported in Supplementary Material.

As a final remark, the use of budgets to link the upper problem with the subproblems might lead to an underestimation of storage capacity. This occurs because the exchange of budgets instead of hourly-resolved flows introduces additional flexibility, which can reduce the system storage requirements. In Section 5.4, we present targeted solutions to address this issue.

*Algorithm 2.* *Budget-based temporal + sectoral/spatial Benders decomposition.*

**Input:**. Set decomposition type (sectoral → $\delta = s$, spatial → $\delta = z$). Set maximum number of iterations $K_{max}$ and convergence tolerance $\varepsilon_{tol}$. Set $y^0 = 0$, $q^0_{w,\delta} = 0$ $\forall w, \delta$.

**Output:** $y^{opt} = 0$, $q^{opt}_{w,\delta} = 0$ $\forall w, \delta$

for $k = 0, \ldots, K_{max}$ do
    for $w \in W$ do
        for $s \in S$ do
            Solve operational subproblem (9) with Eq. (10) if $\delta = s$ or Eq. (11) if $\delta = z$.
        end for
    end for
    Compute best upper bound $UB^k$ as $UB^k = \min_{j=0,\ldots,k} c_y^T y^j + \sum_w \sum_\delta g^j_{w,\delta}$.
    Update cuts in upper problem as in Eq. (12b).
    Solve upper problem (12) to obtain lower bound $LB^k$.
    **if** $(UB^k - LB^k)/LB^k \leq \varepsilon_{tol}$ **then**
        Set $y^{opt} = y^{k+1}$ and $q^{opt}_{w,\delta} = q^{k+1}_{w,\delta}$ $\forall w, \delta$
        **stop**
    **else**
        Solve regularised upper problem (13) to obtain $y^{k+1}$ and $q^{k+1}_{w,\delta}$
    **end if**
end for

## 5. Results

We test the developed BD algorithms against the regular monolithic formulation and the temporal BD, which represents the state of the art for deterministic CEMs. All methods are implemented in the Dolphyn model (He et al., n.d.)[3].

### 5.1. Numerical experiments and computational setup

We consider case studies of the continental United States made up of 16 and 64 zones, defined based on the spatial resolution used in the Integrated Planning Model (IPM) of the Environmental Protection Agency (EPA) (US Environmental Protection Agency, n.d.), as schematised in Figure 5. To test how the different approaches perform varying the temporal resolution, we represent the operational year with hourly resolution using 12, 22, 32, 42, and 52 representative weeks, identified through k-means clustering.

    We consider the electricity and hydrogen sectors, taking data from Shi (2023) and enforcing a net-zero-$CO_2$-emission constraint. In this example application, the connection between the two sectors is limited to power-to-gas flows. Input time series include renewables availability (solar photovoltaic, onshore wind, offshore wind, Run-of-River hydro), fuel prices (natural gas, uranium), electricity demand, and hydrogen demand. The optimisation problem is implemented in Dolphyn, using Julia 1.9.2 (Bezanson et al., 2017) and JuMP 1.20.0 (Dunning et al., 2017). LP problems are solved with Gurobi 10 (Gurobi Optimization, 2014) with the barrier method with crossover disabled, considering a convergence tolerance of $10^{-3}$ for the decomposition algorithm. All numerical experiments are run in the MIT SuperCloud high-performance computing environment (Reuther et al., 2018) with Intel Xeon Platinum 8260 2.40 GHz

---

[3] Data and code will be made available in a public repository upon acceptance of the manuscript.

processors, 48 CPUs per node, and 192 GB RAM per node. For all decomposition algorithms, subproblems are solved in parallel, assigning one CPU to each subproblem. If the number of subproblems exceeds the available CPUs on a single node, distributed parallelisation is used. We compare the runtime of the developed BD algorithms against the regular monolithic formulation and the temporal BD proposed by Jacobson et al. (2024) and Pecci &Jenkins (2024) presented in Section 3, which represents the state of the art for deterministic CEMs. Both the sectoral and spatial BD include decomposition in time, and are therefore referred to as temporal + sectoral and temporal + spatial BD in the remainder of this section. In the method comparisons, runtimes include time spent on both model creation and solution.

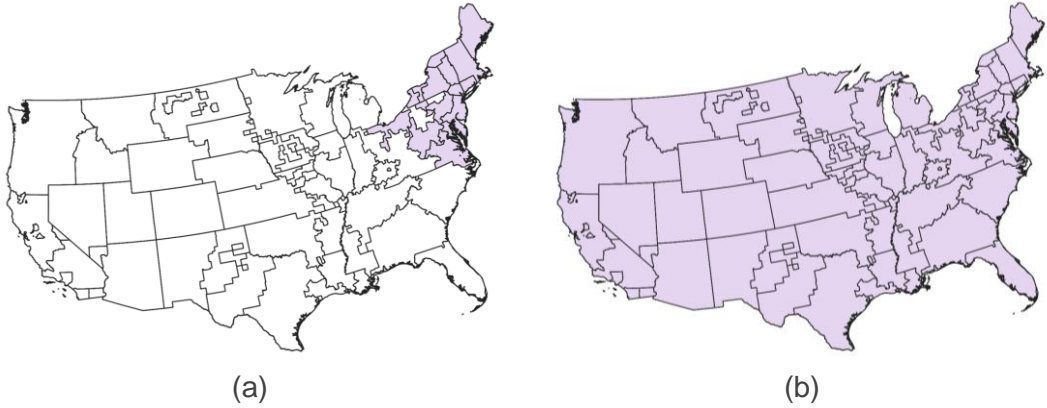

*Figure 5. Spatial configurations considered in the numerical experiments: 16 zones (a), and 64 zones (b).*

5.2.   *Extension of temporal Benders decomposition to multi-sector models*

The temporal BD algorithm is extended to the multi-sector model (electricity and hydrogen) and compared against a single-sector case that considers only electricity. Table 1 presents the computational performance of the temporal BD algorithm and of the monolithic solution. As already observed by Jacobson et al. (2024), decomposition becomes competitive as the model complexity grows. In the single-sector model, temporal BD outperforms the monolithic model in the most complex configuration (64 zones, 52 weeks). The introduction of an additional sector results in a substantial increase in computational complexity, leading to a twentyfold increase in runtime on average. In the 64-zone, 52-week configuration, the runtime of temporal BD increases by a factor of thirteen with the addition of a sector, while the multi-sector monolithic model becomes intractable, running out of memory. This dramatic increase in the computational burden underscores the challenges associated with the application of temporal decomposition to multi-sector models, highlighting the need for novel approaches to improve computational efficiency.

Table 1. Runtime for monolithic models and temporal BD (100 s) for a single-sector (electricity) and a multi-sector (electricity + hydrogen) configuration. The best-performing formulation for each configuration is highlighted in bold. The infinity symbol (∞) is used to indicate intractable problems due to memory.

|  |  | Single-sector | | | | | Multi-sector | | | | |
|---|---|---|---|---|---|---|---|---|---|---|---|
| Representative weeks → | Zones | 12 | 22 | 32 | 42 | 52 | 12 | 22 | 32 | 42 | 52 |
| Monolithic (100 s) | 16 | 1 | 2 | 2 | 3 | 4 | 1 | 5 | 10 | 13 | 20 |
|  | 64 | 2 | 5 | 10 | 13 | 54 | 15 | 59 | 65 | 127 | ∞ |
| Temporal BD (100 s) | 16 | 26 | 4 | 4 | 5 | 6 | 61 | 40 | 39 | 37 | 43 |
|  | 64 | 20 | 36 | 97 | 59 | **46** | 1729 | 1375 | 988 | 934 | **602** |

5.3. Budget-based sectoral Benders decomposition

Table 2 compares the computational performance of the budget-based temporal + sectoral BD against the temporal BD and the monolithic model. Even by further decomposing the problem, implementing decomposition in cases with few representative periods brings little benefits. However, when compared to the temporal BD, the temporal + sectoral BD achieves better performance across all cases, with runtime reductions of 20–70%. For the most complex configuration (64 zones, 52 weeks), where the monolithic model is intractable, the temporal + sectoral BD achieves a 31% runtime reduction compared to the temporal BD.

Table 2. Runtime for monolithic models, temporal BD, and budget-based temporal + sectoral BD (100 s). The best-performing formulation for each configuration is highlighted in bold. The infinity symbol (∞) is used to indicate intractable problems due to memory.

| Representative weeks → | Zones | 12 | 22 | 32 | 42 | 52 |
|---|---|---|---|---|---|---|
| Monolithic (100 s) | 16 | 1 | 5 | 10 | 13 | 20 |
|  | 64 | 15 | 59 | 65 | 127 | ∞ |
| Temporal BD (100 s) | 16 | 61 | 40 | 39 | 37 | 43 |
|  | 64 | 1729 | 1375 | 988 | 934 | 602 |
| Temporal + sectoral BD (100 s) | 16 | 28 | 20 | 18 | 20 | 28 |
|  | 64 | 607 | 376 | 775 | 753 | **416** |

Figure 6 shows the number of iterations and the runtime per iteration for the temporal and temporal + sectoral BD across the considered cases. By doubling the cuts introduced at each iteration, the temporal + sectoral BD achieves convergence in fewer iterations. The runtime per iteration does not exhibit a clear trend due to the interplay of two opposing effects. On one hand, the temporal + sectoral BD further decomposes the problem, resulting in smaller subproblems that can be solved in parallel, potentially reducing the computational time per iteration. On the other hand, the increased number of cuts leads to a larger upper problem, affecting its solution time. Results shows that, in the considered case studies, the reduction in subproblem size dominates when the number of representative periods is smaller, leading to consistently lower runtime per iteration for the temporal + sectoral BD. This is particularly evident in the 64-zone configuration with 12 and 22 representative weeks, where runtime per iteration is 50–70% lower than that of the temporal BD. In contrast, for configurations with a high number of representative periods, the two opposing effects balance out, resulting in similar runtimes per iteration. As a result, while the temporal + sectoral BD outperforms the temporal BD across all cases, the runtime reduction is more pronounced as the number of representative periods decreases.

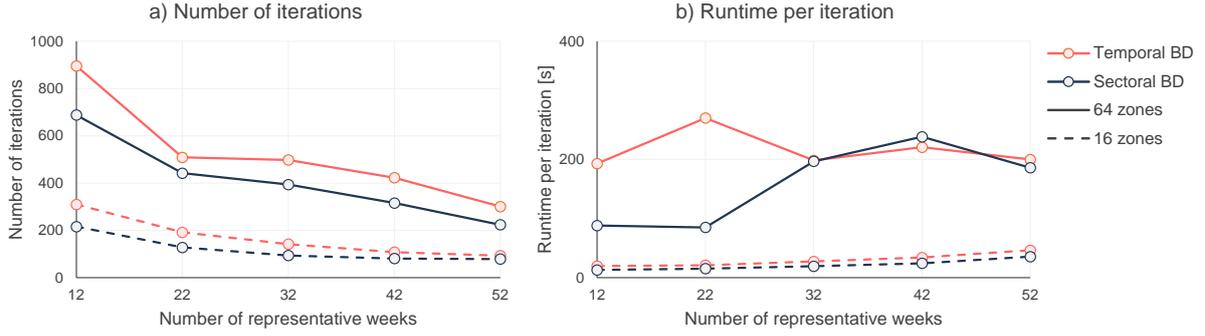

Figure 6. Number of iterations and runtime per iteration for temporal and temporal + sectoral BD in the analysed configurations.

Figure 7 illustrates the accuracy of the two decomposition algorithms in terms of relative error compared to the solution of the monolithic model of the corresponding configuration. Results indicate that the budget-based temporal + sectoral BD achieves accuracy comparable to the temporal BD, often providing better estimates for hydrogen production. The limitation of the budget-based temporal + sectoral BD lies in the estimation of the hydrogen storage capacity, which results consistently equal to zero. This issue is inherent to the budget-based formulation. By exchanging budgets, the hourly export profile of a vector from one sector may not match with its hourly import profile in another. In the example application, this results in added flexibility in electricity consumption for electrolysis, which reduces the need for hydrogen storage. Consequently, additional steps are required to estimate hydrogen storage capacity. This issue is investigated in the next section, which proposes a solution to address it.

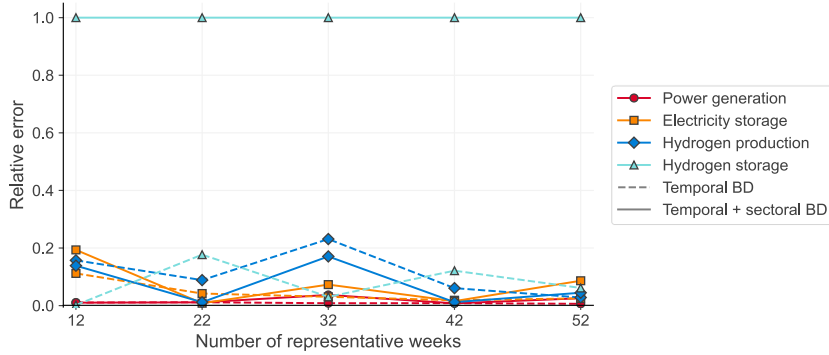

Figure 7. Relative error of installed capacities for temporal (dashed lines) and temporal + sectoral (solid lines) BD compared to solutions of the monolithic formulation. Results are shown for the 16-zone configuration.

### 5.4. Two-stage sectoral Benders decomposition

The developed budget-based sectoral BD achieves high accuracy for capacity expansion of all system components, but it fails to estimate hydrogen storage requirements. To address this limitation, we introduce a two-stage algorithm. In the first stage, we solve the budget-based temporal + sectoral BD (Algorithm 2), and we use the computed cuts to warm start the sector-aggregated problem of Algorithm 1, connecting the θ variables of the two upper problems as:

$$\forall w \quad \theta_w \geq \sum_{s \in S} \theta_{w,s} \qquad (14)$$

To manage the potentially large number of cuts in the final upper problem of the budget-based decomposition, non-binding cuts are removed. These are identified as the ones with a large slack value ($10^1$ in this application), computed by evaluating the difference between the left-hand side and right-hand side of each cut after solving the problem. Given the accuracy of the budget-based decomposition, the obtained capacities of power generation, electricity storage, and hydrogen production units are used to define bounds for the corresponding decision variables in the second stage of the algorithm. Specifically, the resulting capacities are set as lower bound, while an upper bound is set with a +5% margin. Since the second stage is exclusively aimed at estimating hydrogen storage capacity (while capacities of all other system components are accurately determined in the first stage), a coarser convergence tolerance can be used. In the considered applications, we use a tolerance of $10^{-2}$, which has proven sufficient to achieve good accuracy on the storage capacity. The resulting algorithm is summarised in Algorithm 3.

*Algorithm 3.* Two-stage temporal + sectoral Benders decomposition.

**Stage 1:** Solve Algorithm 2.

**Stage 2:** Solve Algorithm 1 initialising Problem (4) with the cuts computed in Stage 1, adding Eq. (14), and removing all non-binding cuts.

Runtimes for the two-stage temporal + sectoral BD are reported in Table 3 and compared to the temporal BD. Despite the additional computational stage, the temporal + sectoral BD continues to outperform the temporal BD across all configurations, with runtime reductions within 15-70%. On average, the second stage requires 10-20 iterations to achieve convergence. As Figure 8 shows, the introduction of the second stage effectively addresses the limitations of the budget-based formulation, enabling a good estimation of hydrogen storage requirements. As a result, the accuracy of the temporal + sectoral BD becomes comparable to that of the temporal BD, with relative errors remaining below 20% across all cases. This level of accuracy is considered acceptable, given that multiple near-optimal solutions with different configurations but similar total annual costs may exist (DeCarolis et al., 2016; Pickering et al., 2022).

Table 3. Runtime for temporal BD and two-stage temporal + sectoral BD (100 s). The best-performing formulation for each configuration is highlighted in bold.

| Representative weeks → | Zones | 12 | 22 | 32 | 42 | 52 |
|---|---|---|---|---|---|---|
| Temporal BD (100 s) | 16 | 61 | 40 | 39 | 37 | 43 |
| | 64 | 1729 | 1375 | 988 | 934 | 602 |
| Two-stage temporal + sectoral BD (100 s) | 16 | **49** | **26** | **22** | **24** | **30** |
| | 64 | **501** | **406** | **817** | **808** | **510** |

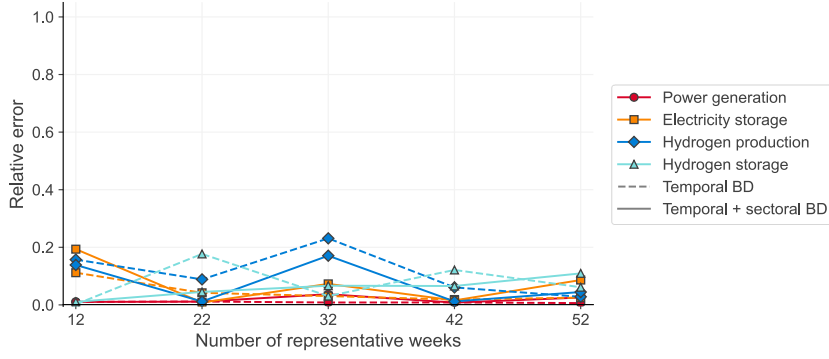

*Figure 8. Relative error of installed capacities for temporal (dashed lines) and two-stage temporal + sectoral (solid lines) BD compared to solutions of the monolithic formulation. Results are shown for the 16-zone configuration.*

### 5.5. Budget-based spatial Benders decomposition

The decomposition of the problem in both time and space results in a large number of subproblems. While this increases the number of cuts, improving convergence, it also introduces significant complexity in managing distributed computing, often leading to out-of-memory errors. For this reason, we here test the method on a single case study as a proof of concept to demonstrate its computational efficiency and potential for significant runtime reduction, while the development of additional strategies to address these challenges is left for future work.

Figure 9 compares the convergence of the temporal + spatial, temporal + sectoral, and temporal BD for the 16-zone, 12-week configuration. Algorithm 2 is modified by removing non-binding cuts every 10 iterations in order to limit the size of the upper problem, which might increase rapidly due to addition of numerous cuts per iteration. The decomposition of the problem in space enables a more than tenfold reduction in the number of iterations, from 310 and 216 for the temporal and temporal + sectoral BD, respectively, to only 20 for the temporal + spatial BD. This translates into a runtime reduction of approximately 70% compared to the temporal BD and 40% compared to the temporal + sectoral BD.

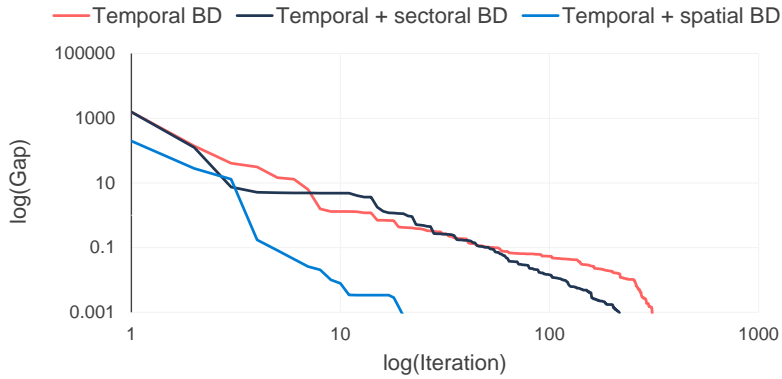

*Figure 9. Convergence of temporal + spatial, temporal + sectoral, and temporal BD for the 16-zone, 12-week configuration.*

As in the temporal + sectoral BD, budget-based linking results in an underestimation of the storage requirements. To address this issue, we adopt the approach presented in Section 5.4, implementing a two-stage algorithm. In the first stage, we solve the budget-based temporal + spatial BD (Algorithm 2), and we use the computed cuts to warm start the spatially aggregated problem of Algorithm 1. The θ variables of the two upper problems are connected as:

$$\forall w \qquad \theta_w \geq \sum_{z \in Z} \theta_{w,z} \qquad (15)$$

The resulting algorithm is summarised in Algorithm 4. Similar to the two-stage temporal + sectoral BD, the second stage primarily aims to estimate storage capacities. Therefore, a coarser tolerance ($10^{-2}$ in the considered applications) can be used.

---

***Algorithm 4.*** *Two-stage temporal + spatial Benders decomposition.*

---

**Stage 1:** Solve Algorithm 2.

**Stage 2:** Solve Algorithm 1 initialising Problem (4) with the cuts computed in Stage 1, adding Eq. (15), and removing all non-binding cuts.

---

In the investigated configuration, the second stage requires an additional 21 iterations. This leads to a total of 41 iterations, which remains one order of magnitude lower than the temporal + sectoral and temporal BD, as schematised in Figure 10.a. Similarly, the runtime increases only marginally, achieving a runtime reduction of approximately 50–60% compared to the two-stage temporal + sectoral BD and temporal BD (Figure 10.b).

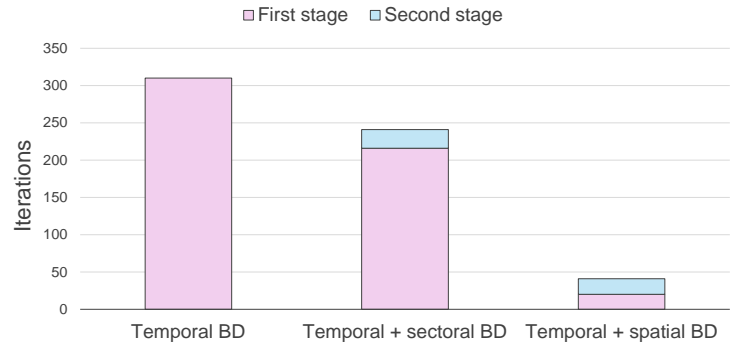

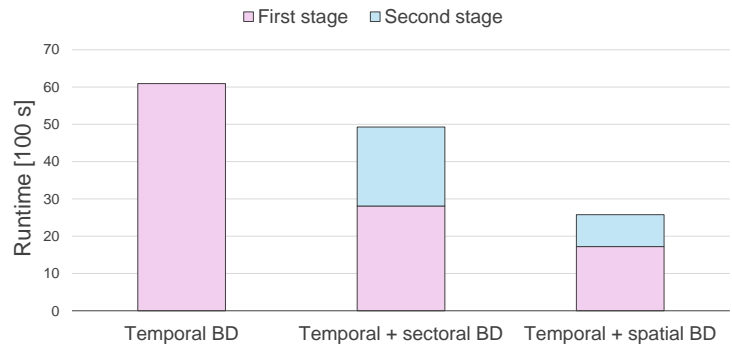

*Figure 10. Number of iterations (a) and runtime (b) for temporal, two-stage temporal + sectoral, and two-stage temporal + spatial BD in the 16-zone, 12-week configuration. Since the temporal BD does not require a second stage, all iterations and the entire runtime are allocated to the first stage.*

Figure 11 shows the accuracy of the approach by comparing installed capacities and the objective function value normalised to the monolithic model solution. The introduction of the second stage allows for the estimation of electricity and hydrogen storage capacity, which would otherwise be absent.

Compared to the monolithic model, results of the two-stage temporal + spatial BD feature higher power generation and electricity storage capacity, while hydrogen production and hydrogen storage capacities are lower. However, the objective function deviates by only 0.5%, indicating that two system configurations achieve comparable total annual cost through different pathways, with the temporal + spatial BD solution being more balanced towards the electricity sector. This outcome aligns with findings in the energy system modelling literature, which suggest that multiple system configurations can achieve similar total costs (DeCarolis et al., 2016; Pickering et al., 2022). Accordingly, the accuracy of the proposed approach is considered acceptable.

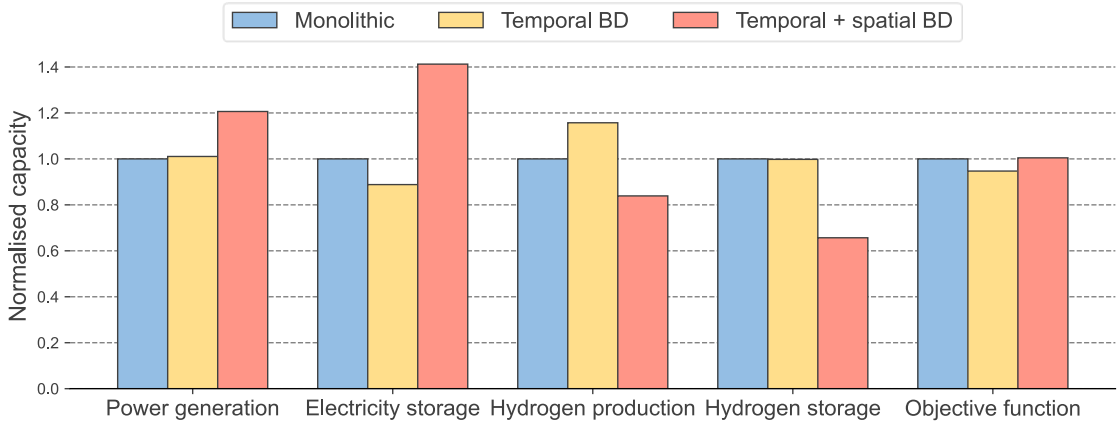

Figure 11. Cumulative installed capacities (power generation, electricity storage, hydrogen production, and hydrogen storage) and objective function value normalised to the monolithic model solution, in the 16-zone, 12-week configuration.

## 6. Conclusions

In this work, we developed sectoral and spatial Benders decomposition algorithms to improve the computational performance of large-scale capacity expansion models used for energy planning. The key novelty lies in the budget-based formulation, which leverages the structure of this class of models to enable efficient linking between the upper problem and subproblems, improving the convergence of the Benders algorithm. The budget-based Benders decomposition algorithms demonstrated substantial improvements of the computational performance over state-of-the-art decomposition methods. We considered as reference a state-of-the-art temporal Benders decomposition algorithm and extended the decomposition to the sectoral or spatial domain. The temporal + sectoral Benders decomposition achieves 20-70% runtime reductions with respect to the temporal decomposition for a system with electricity and hydrogen as sectors. The approach shows potential for further efficiency improvement when extended to additional sectors, providing additional cuts per iteration. The temporal + spatial Benders decomposition achieves even greater performance gains. In the tested application, the number of iterations is reduced by an order of magnitude with respect to both the temporal and temporal + sectoral decomposition, leading to a 40-70% reduction in runtime. While spatial decomposition introduces additional challenges in managing a large number of subproblems in a distributed computing environment, the substantial gain in computational performance demonstrates its potential, fully justifying further efforts in developing strategies to address this limitation.

Relying on a budget-based formulation leads to an underestimation of storage requirements. To overcome this limitation, we introduced a two-stage algorithm, which builds on the solution of the budget-based decomposition to incorporate an additional computational step for storage capacity estimation. We applied the two-stage algorithm to both the sectoral and the spatial Benders decomposition, showing that the approach effectively addresses the issue while preserving superior computational efficiency, with runtime reductions compared to the temporal Benders decomposition remaining in the range of 15–70%.

To fully exploit the advantages of these decomposition algorithms, future developments should focus on strategies to efficiently manage distributed computing with a large number of subproblems. A possible approach is to apply spatial decomposition at the level of zone clusters rather than individual zones. While this would reduce some of the advantages of spatial decomposition, it would also alleviate the computational burden by limiting the number of subproblems that need to be solved in parallel. Overall, the analysis demonstrates the potential of the developed budget-based decomposition formulations - both sectoral and spatial - to significantly improve the computational efficiency of capacity expansion models, enabling model tractability with higher temporal, spatial, and technological resolution. While in this work we tested temporal decomposition combined with either sectoral or spatial decomposition, the approaches can be integrated into a unified temporal + sectoral + spatial decomposition algorithm, combining their advantages to further improve the computational performance. With respect to generalisation, the proposed methods leverage the generic structure of capacity expansion models, making them broadly applicable to most existing integrated energy system models.

## Acknowledgements

The authors thankfully acknowledge the Progetto Rocca program for supporting this collaboration between MIT and Politecnico di Milano. F. Parolin gratefully acknowledges the MIT Energy Initiative for their hospitality and support.

# Supplementary Material for

# Sectoral and spatial decomposition methods for multi-sector capacity expansion models

*Federico Parolin [a,b], Yu Weng [b], Paolo Colbertaldo [a], Ruaridh Macdonald [b]*

a. Department of Energy, Politecnico di Milano, Milan, Italy
b. MIT Energy Initiative, Massachusetts Institute of Technology, Cambridge, MA, USA

## 1. Formulation: Budget-based sectoral Benders decomposition

The balance of an energy vector $v$ in sector $s$ at zone $z$ and time step $t$ can be expressed as:

$$\sum_g x_{gen}^{g,v,z,s,t} + \sum_\sigma \frac{x_{otp}^{\sigma,v,z,s,t}}{\eta_{\sigma,otp}} + x_{nse}^{v,z,s,t} - x_{dem}^{v,z,s,t} - \sum_{z'\neq z} x_{trn}^{v,z,z',t} \sum_\sigma x_{ipt}^{\sigma,v,z,s,t} \cdot \eta_{\sigma,ipt} - x_{crt}^{v,z,s,t} = \sum_{s'\neq s} x_{exp}^{v,z,s,s',t} \quad (1)$$

where $x_{gen}^{g,v,z,s,t}$ is the energy vector generation from technology $g$, $x_{otp}^{\sigma,v,z,s,t}$ and $x_{ipt}^{\sigma,v,z,s,t}$ are the output and input flow from storage technology $\sigma$, $x_{trn}^{v,z,z',t}$ is the transport flow from zone $z$ to zone $z'$, $x_{nse}^{v,z,s,t}$ is the non-served energy, $x_{dem}^{v,z,s,t}$ is the energy vector demand, $x_{crt}^{v,z,s,t}$ is the curtailment, and $x_{exp}^{v,z,s,s',t}$ is the net energy vector export from sector $s$ to sector $s'$.

This structure allows to highlight the role of net export flows, which serve as the connecting variables between different sectors. These enter an additional balance, which ensures that the net export from sector $s$ and sector $s'$ is equal and opposite to the net export from sector $s'$ to sector $s$:

$$\forall v,z,t,s\neq s' \qquad x_{exp}^{v,z,s,s',t} = -x_{exp}^{v,z,s',s,t} \quad (2)$$

We introduce export budgets for every subperiod as complicating variables, defined as:

$$\forall z, s\neq s' \qquad y_{exp}^{v,z,s,s',w} = \sum_{t\in T_w} x_{exp}^{v,z,s,s',t} \quad (3)$$

Export budget are treated as investment decision variables and included in the set $y \in Y$. We define the operational subproblem as:

$$g_{w,s}^k = \min \sum_{z \in Z} c_{w,z,s}^T x_{w,z,s} \tag{4a}$$

$$\text{s.t.} \quad A_{w,z,s} x_{w,z,s} + B_{w,z,s} y \leq b_{w,s,z} \quad \forall z \in Z \tag{4b}$$

$$\sum_{z \in Z} Q_{w,z,s} x_{w,z,s} \leq q_{w,s} \tag{4c}$$

$$\sum_{t \in T_w} x_{\exp}^{v,z,s,s',t} = y_{\exp}^{v,z,s',s,w} \quad \forall z \in Z \tag{4d}$$

$$y = y^k \quad : \pi^k \tag{4e}$$

$$q_{w,s} = q_w^k \quad : \lambda_{w,s}^k \tag{4f}$$

$$x_{w,z,s} \geq 0 \quad \forall z \in Z \tag{4g}$$

which is solved for each subperiod and sector.

The upper problem is defined as:

$$LB^k = \min c_y^T y^j + \sum_{w \in W} \sum_{s \in S} \theta_{w,s} \tag{5a}$$

$$\text{s.t.} \quad \theta_{w,s} \geq g_{w,s}^j + (y - y^j)^T \pi^j + \left(q_{w,s} - q_{w,s}^j\right)^T \lambda_{w,s}^j \quad \forall j = 0, \dots, k, w \in W, s \in S \tag{5b}$$

$$\sum_{w \in W} \sum_{s \in S} q_{w,s} \leq e \tag{5c}$$

$$y_{\exp}^{v,z,s,s',w} = -x_{\exp}^{v,z,s',s,w} \quad \forall w \in W, z \in Z, s \neq s' \tag{5d}$$

$$y_{\exp}^{v,z,s,s',w} \leq \sum_{g \in G} y_{\text{gen}}^{g,v,s,z} \cdot |T_w| \quad \forall w \in W, z \in Z, s \neq s' \tag{5e}$$

$$Ry \leq r \tag{5f}$$

$$y \in Y \tag{5g}$$

where Eq. (5e) is introduced to avoid large values of budgets in the subproblems. The regularisation problem is updated accordingly, as:

$$\min \Phi^{\text{int}} = 0 \tag{6a}$$

$$\text{s.t.} \quad \theta_{w,s} \geq g_{w,s}^j + (y - y^j)^T \pi^j + \left(q_{w,s} - q_{w,s}^j\right)^T \lambda_{w,s}^j \quad \forall j = 0, \dots, k, w \in W, s \in S \tag{6b}$$

$$\sum_{w \in W} \sum_{s \in S} q_{w,s} \leq e \tag{6c}$$

$$y_{\exp}^{v,z,s,s',w} = -y_{\exp}^{v,z,s',s,w} \quad \forall w \in W, z \in Z, v \in V, s \neq s' \tag{6d}$$

$$y_{\exp}^{v,z,s,s',w} \leq \sum_{g \in G} y_{\text{gen}}^{g,v,s,z} \quad \forall w \in W, z \in Z, s \neq s' \tag{6e}$$

$$Ry \leq r \tag{6f}$$

$$y \in Y \tag{6g}$$

$$c_y^T y^j + \sum_{w \in W} \sum_{s \in S} \theta_{w,s} \leq LB^k + \alpha(UB^k - LB^k) \tag{6h}$$

The full solution algorithm is reported in Algorithm 1.

---

*Algorithm 1.* Sectoral Benders decomposition.

---
**Input:** $y^0 = 0$, $q_{w,s}^0 = 0$ $\forall w \in W, s \in S$. Set maximum number of iterations $K_{\max}$ and convergence tolerance $\varepsilon_{\text{tol}}$.

**Output:** $y^{\text{opt}} = 0$, $q_{w,s}^{\text{opt}} = 0$ $\forall w \in W, s \in S$

**for** $k = 0, \dots, K_{\max}$ **do**

    **for** $w \in W$ **do**

        **for** $s \in S$ **do**

            Solve operational subproblem (4).

        **end for**

    **end for**

    Compute best upper bound $UB^k$ as $UB^k = \min_{j=0,\dots,k} c_y^T y^j + \sum_{w \in W} \sum_{s \in S} g_{w,s}^j$.

    Update cuts in upper problem as in Eq. (5.b).

    Solve upper problem (5) to obtain lower bound $LB^k$.

    **if** $(UB^k - LB^k)/LB^k \leq \varepsilon_{\text{tol}}$ **then**

        Set $y^{\text{opt}} = y^{k+1}$ and $q_{w,s}^{\text{opt}} = q_{w,s}^{k+1}$ $\forall w \in W, s \in S$

        **stop**

    **else**

        Solve regularised upper problem (6) to obtain $y^{k+1}$ and $q_{w,s}^{k+1}$

    **end if**

**end for**

---

## 2. Formulation: Budget-based spatial Benders decomposition

The zonal balance of Eq. (1) can be divided to isolate the contribution of transport flows, as:

$$\sum_g x_{\text{gen}}^{g,v,z,s,t} + \sum_\sigma \frac{x_{\text{otp}}^{\sigma,v,z,s,t}}{\eta_{\sigma,\text{otp}}} + x_{\text{nse}}^{v,z,s,t} - x_{\text{dem}}^{v,z,s,t} - \sum_\sigma x_{\text{ipt}}^{\sigma,v,z,s,t} \cdot \eta_{\sigma,\text{ipt}} - x_{\text{crt}}^{v,z,s,t} - \sum_{s' \neq s} x_{\exp}^{v,z,s,s',t} = \sum_{z' \neq z} x_{\text{trn}}^{v,z,z',t} \tag{7}$$

We introduce budgets for transport flows and we use them as complicating variables in the decomposition, defining them as:

$$\forall z, v, w \qquad y_{trn}^{v,z,w} = \sum_{t \in T_w} \sum_{z' \neq z} x_{trn}^{v,z,z',t} \tag{8}$$

We define spatial subproblems as:

$$g_{w,z}^k = \min \sum_{s \in S} c_{w,z,s}^T x_{w,z,s} \tag{9a}$$

$$\text{s.t.} \quad A_{w,z,s} x_{w,z,s} + B_{w,z,s} y \leq b_{w,s,z} \quad \forall s \in S \tag{9b}$$

$$\sum_{s \in S} Q_{w,z,s} x_{w,z,s} \leq q_{w,z} \tag{9c}$$

$$\sum_{t \in T_w} \sum_{z' \neq z} x_{trn}^{v,z,z',t} = y_{trn}^{v,z,w} \quad \forall v \in V, w \in W \tag{9d}$$

$$y = y^k \quad : \pi^k \tag{9e}$$

$$q_{w,z} = q_w^{\ k} \quad : \lambda_{w,z}^k \tag{9f}$$

$$x_{w,z,s} \geq 0 \quad \forall s \in S \tag{9g}$$

which are solved for each subperiod and zone. The new upper problem is built as:

$$LB^k = \min c_y^T y^j + \sum_{w \in W} \sum_{z \in Z} \theta_{w,z} \tag{10a}$$

$$\text{s.t.} \quad \theta_{w,z} \geq g_{w,z}^j + \left(y - y^j\right)^T \pi^j + \left(q_{w,z} - q_{w,z}^j\right)^T \lambda_{w,z}^j \quad \forall j = 0, \ldots, k, w \in W, z \in z \tag{10b}$$

$$\sum_{w \in W} \sum_{z \in Z} q_{w,z} \leq e \tag{10c}$$

$$\sum_{t \in T_w} y_{trn}^{v,z,t} = y_{trn}^{v,z,w} \quad \forall v \in V, w \in W \tag{10d}$$

$$y_{trn}^{v,z,t} = \sum_{z' \neq z} x_{trn}^{v,z,z',s,t} \quad \forall z \in Z, s \in S, v \in V, t \in T_{w \in W} \tag{10e}$$

$$Ry \leq r \tag{10f}$$

$$y \in Y \tag{10g}$$

where we introduce a dedicated variable $y_{trn}^{v,z,t}$ to translate budgets into hourly flows, enabling the definition of the transport flow balance in the upper problem (Eq. (10e)). The regularisation problem is updated accordingly, as:

$$\min \Phi^{int} = 0 \tag{11a}$$

$$\text{s.t.} \quad \theta_{w,z} \geq g^j_{w,z} + (y - y^j)^T \pi^j + (q_{w,z} - q^j_{w,z})^T \lambda^j_{w,z} \quad \forall j = 0, \ldots, k, w \in W, z \in z \tag{11b}$$

$$\sum_{w \in W} \sum_{z \in Z} q_{w,z} \leq e \tag{11c}$$

$$\sum_{t \in T_w} y^{v,z,t}_{trn} = y^{v,z,w}_{trn} \quad \forall v \in V, w \in W \tag{11d}$$

$$y^{v,z,t}_{trn} = \sum_{z' \neq z} x^{v,z,z',t}_{trn} \quad \forall v \in V, t \in T_{w \in W} \tag{11e}$$

$$Ry \leq r \tag{11f}$$

$$y \in Y \tag{11g}$$

$$c_y^T y^j + \sum_{w \in W} \sum_{z \in Z} \theta_{w,z} \leq LB^k + \alpha(UB^k - LB^k) \tag{11h}$$

The resulting solution algorithm is reported in Algorithm 2.

---

**Algorithm 2.** *Spatial Benders decomposition.*

---

**Input:** $y^0 = 0, q^0_{w,z} = 0 \ \forall w \in W, z \in Z$. Set maximum number of iterations $K_{max}$ and convergence tolerance $\varepsilon_{tol}$.

**Output:** $y^{opt} = 0, q^{opt}_{w,z} = 0 \ \forall w \in W, z \in Z$

**for** $k = 0, \ldots, K_{max}$ **do**

    **for** $w \in W$ **do**

        **for** $z \in Z$ **do**

            Solve operational subproblem (9).

        **end for**

    **end for**

    Compute best upper bound $UB^k$ as $UB^k = \min_{j=0,\ldots,k} c_y^T y^j + \sum_{w \in W} \sum_{z \in Z} g^j_{w,z}$.

    Update cuts in upper problem as in Eq. (10b).

    Solve upper problem (10) to obtain lower bound $LB^k$.

    **if** $(UB^k - LB^k)/LB^k \leq \varepsilon_{tol}$ **then**

        Set $y^{opt} = y^{k+1}$ and $q^{opt}_{w,z} = q^{k+1}_{w,z} \ \forall w \in W, z \in Z$

        **stop**

    **else**

        Solve regularised upper problem (11) to obtain $y^{k+1}$ and $q^{k+1}_{w,z}$

    **end if**

**end for**

---